\documentclass[12pt]{article}
\usepackage{amsmath}
\usepackage{amssymb}
\usepackage[dvips]{graphicx}
\usepackage{epsfig}
\usepackage{subfigure}
\usepackage{vector}
\def\eps{\varepsilon}
\font\tencmmib=cmmib10 \skewchar\tencmmib '60
\newfam\cmmibfam
\textfont\cmmibfam=\tencmmib

\def\bbox{\quad\hbox{\vrule \vbox{\hrule \vskip2pt \hbox{\hskip2pt
\vbox{\hsize=1pt}\hskip2pt} \vskip2pt\hrule}\vrule}}
\def\lessim{\ \lower4pt\hbox{$
\buildrel{\displaystyle <}\over\sim$}\ }
\def\gessim{\ \lower4pt\hbox{$\buildrel{\displaystyle >}
\over\sim$}\ }

\def\P{{\cal P}}

\def\I{{\rm I}}
\def\II{{\rm II}}

\def\M{{\cal M}}

\def\C{{\cal C}}
\def\Cp{{{\cal C}^{\prime}}}
\def\eps{{\varepsilon}}
\def\ch{{\mbox{\rm ch}}}

\def\qed{\hfill\break\rightline{$\bbox$}}
\newcommand{\e}{\mathbb{E}}
\newcommand{\p}{\mathbb{P}}
\newcommand{\Reals}{\mathbb{R}}

\newtheorem{lemma}{Lemma}
\newtheorem{theorem}{Theorem}
\newtheorem{corollary}{Corollary}
\makeatletter
\@addtoreset{equation}{section}

\makeatother

\font\tencmmib=cmmib10 \skewchar\tencmmib '60
\newfam\cmmibfam
\textfont\cmmibfam=\tencmmib

\def\bbox{\quad\hbox{\vrule \vbox{\hrule \vskip2pt \hbox{\hskip2pt
\vbox{\hsize=1pt}\hskip2pt} \vskip2pt\hrule}\vrule}}
\def\lessim{\ \lower4pt\hbox{$
\buildrel{\displaystyle <}\over\sim$}\ }
\def\gessim{\ \lower4pt\hbox{$\buildrel{\displaystyle >}
\over\sim$}\ }

\def\eps{\varepsilon}

\def\go0{\to 0}
\def\grad{\triangledown}

\def\leftitem#1{\item{\hbox to\parindent{\enspace#1\hfill}}}

\def\qed{\hfill\break\rightline{$\bbox$}}

\def\sg{\sigma}

\def\sg2{\sigma^2}

\def\__{_{\infty}}
\begin{document}

\title{
A question about Parisi functional.}

\author{ 
Dmitry Panchenko\thanks{Department of Mathematics, Massachusetts Institute
of Technology, 77 Massachusetts Ave, Cambridge, MA 02139
email: panchenk@math.mit.edu}\\
{\it Department of Mathematics}\\
{\it Massachusetts Institute of Technology}\\
}

\maketitle
\begin{abstract}
We conjecture that the Parisi functional in the Sherrington-Kirkpatrick
model is convex in the functional order parameter. We prove a partial result
that shows the convexity along ``one-sided'' directions. 
An interesting consequence of this result is the log-convexity of
$L_m$ norm for a class of random variables.

\end{abstract}
\vspace{0.5cm}


\section{A problem and some results.}

Let $\M$ be a set of all nondecreasing and right-continuous
functions $m:[0,1]\to [0,1].$
Let us consider two convex smooth functions
$\Phi$ and $\xi:\Reals\to \Reals$ both symmetric, 
$\Phi(-x) = \Phi(x)$ and $\xi(-x)=\xi(x),$
and $\Phi(0)=\xi(0)=0.$ 
%Define $\theta(x)=x\xi'(x)-\xi(x).$
We will also assume that $\Phi$ is of moderate growth so that
all integrals below are well defined.

Given $m\in\M,$ consider a function $\Phi(q,x)$ for
$q\in [0,1], x\in\Reals$ such that $\Phi(1,x) = \Phi(x)$ and
\begin{equation}
\frac{\partial \Phi}{\partial q} = -\frac{1}{2}\xi''(q)
\Bigl(
\frac{\partial^2 \Phi}{\partial x^2} + m(q)
\Bigl(\frac{\partial \Phi}{\partial x}\Bigr)^2
\Bigr).
\label{heat}
\end{equation}
Let us consider a functional $\P:\M\to \Reals$ defined by
$\P(m) = \Phi(0,h)$ for some $h\in\Reals.$

{\bf Main question:} Is $\P$ a convex functional on $\M$?

The same question was asked in \cite{T-Par}. Unfortunately, 
despite considerable effort, we were not able to give complete 
answer to this question.
In this note we will present a partial result that shows convexity along
the directions $\lambda m + (1-\lambda) n$ when $m(q)\geq n(q)$
for all $q\in [0,1].$ It is possible that the answer to this question
lies in some general principle that we are not aware of. 
A good starting point would be to find an alternative proof of
the simplest case of constant $m$ given in Corollary 1 below.

The functional $\P$ arises in the Sherrington-Kirkpatrick mean field model 
where with the choice of $\Phi(x) = \log\ch x,$
the following {\it Parisi formula}
\begin{equation}
\inf_{m\in\M}\Bigl(\log 2 + \P(m)-\frac{1}{2}\int_{0}^{1} m(q)q\xi''(q)dq\Bigr)
\label{pari}
\end{equation}
gives the free energy of the model. A rigorous proof of this result was
given by Michel Talagrand in \cite{T-P}.
Since the last term is a linear  
functional of $m,$ convexity of $\P(m)$ would imply the uniqueness
of {\it the functional order parameter} $m(q)$ that minimizes 
(\ref{pari}). A particular case of $\xi(x) = \beta^2 x^2/2$ 
for $\beta>0$ would also be of interest since it corresponds to
the original SK model \cite{SherK}.

In the case when $m$ is a step function, the solution of (\ref{heat})
can be written explicitly, since for a constant $m$ the function
$g(q,x) = \exp m \Phi(q,x)$ satisfies the heat equation
$$
\frac{\partial g}{\partial q} = -\frac{1}{2}\xi''(q) 
\frac{\partial^2 g}{\partial x^2}.
$$
Given $k\geq 1,$ let us consider a sequence
$$
0=m_0 \leq m_1\leq\ldots \leq m_k=1
$$
and a sequence
$$
q_0=0 \leq q_1\leq \ldots \leq q_k\leq q_{k+1}=1.
$$
We will denote $\vec{m}=(m_0,\ldots,m_k)$ and $\vec{q}=(q_0,\ldots,q_{k+1}).$
Let us define a function $m\in\M$ by
\begin{equation}
m(q)=m_l \mbox{ for } q_l \leq q < q_{l+1}.
\label{step}
\end{equation}
For this step function $\P(m)$ can be defined as follows.
Let us consider a sequence of independent
Gaussian random variables $(z_l)_{0\leq l\leq k}$ such that
$$
\e z_l^2 = \xi'(q_{l+1}) - \xi'(q_l).
$$ 
Define $\Phi_{k+1}(x) = \Phi(x)$ and recursively
over $l\geq 0$ define
\begin{equation}
\Phi_l(x) = \frac{1}{m_l}\log \e_l \exp m_l \Phi_{l+1}(x+z_l)
\label{rec}
\end{equation}
where $\e_l$ denotes the expectation in $(z_i)_{i\geq l}$
and in the case of $m_l=0$ this means $\Phi_l(x) = 
\e_l \Phi_{l+1}(x+z_l).$
Then $\P(m)$ for $m$ in (\ref{step}) is be given by
\begin{equation}
\P_k = \P_k(\vec{m},\vec{q}) = \Phi_0(h).
\label{Pk}
\end{equation}
For simplicity of notations, we will sometimes omit the dependence
of $\P_k$ on $\vec{q}$ and simply write $\P_k(\vec{m}).$
Let us consider another sequence $\vec{n} = (n_0,\ldots,n_k)$ such that
$$
0=n_0\leq n_1\leq \ldots \leq n_k =1.
$$
The following is our main result.

\begin{theorem}\label{T2}
If $n_j \leq m_j$ for all $j$ 
or $n_j \geq m_j$ for all $j$ then
\begin{equation}
\P_k(\vec{n}) - \P_k(\vec{m}) \geq 
\grad \P_k(\vec{m})\cdot(\vec{n} - \vec{m}) =
\sum_{0\leq j\leq k}\frac{\partial \P_k}{\partial m_j}(\vec{m})(n_j-m_j).
\label{grad}
\end{equation}
\end{theorem}

Since the functional $\P$ is uniformly continuous on $\M$ with respect
to $L_1$ norm (see \cite{Guerra} or \cite{T-Par}), 
approximating any function 
by the step functions implies that $\P$ is continuous along 
the directions $\lambda m + (1-\lambda) n$ when $m(q)\geq n(q)$
for all $q\in [0,1].$

Of course, (\ref{grad}) implies that $\P_k(\vec{m})$ is convex
in each coordinate. This yields an interesting consequence
for the simplest case of a constant function $m(q)=m,$
which formally corresponds to the case of $k=2,$
$$
0=m_0\leq m \leq m_2=1 \mbox{ and } 
0=q_0=q_1 \leq q_2=q_3=1.
$$
In this case,
\begin{equation}
\P_k= f(m) = \frac{1}{m}\log \e \exp m\Phi(h+\sigma z).
\label{onem}
\end{equation}
Here $\sigma^2 =\xi'(1)$ can be made arbitrary by the choice of $\xi.$
(\ref{grad}) implies the following.

\begin{corollary} 
If $\Phi(x)$ is convex and symmetric them $f(m)$ defined in
(\ref{onem}) is convex.
\end{corollary}

{\bf Remark.} In Corollary 1 we do not restrict $m$ to be in $[0,1].$
Also, Theorem \ref{T2} holds without the restriction $m_k=n_k=1.$
We made this assumption only because our main interest is in 
(\ref{pari}) where this restriction is assumed.

Corollary 1 implies that the $L_m$ norm of $\exp \Phi(h+\sigma z)$
is log-convex in $m.$ This is a stronger statement than the well-known
consequence of H\"older's inequality that the $L_m$ norm is
always log-convex in $1/m.$
At this point it does not seem obvious how to give an easier proof
even in the simplest case of Corollary 1 than the one we give below.
For example, it is not clear how to show directly that 
$$
f''(m)= m^{-3}(\e V\log^2 V - (\e V\log V)^2 -2\e V\log V)\geq 0,
$$
where $V=\exp m(\Phi(h+\sigma z) - f(m)).$ 

Finally, let us note some interesting consequences of the convexity
of $f(m).$
First, $f''(0)\geq 0$ implies that the third cumulant of 
$\eta = \Phi(h+\sigma z)$ is nonnegative,
\begin{equation}
\e \eta^3 - 3\e\eta^2 \e\eta
+2(\e\eta)^3 \geq 0.
\label{cumulant}
\end{equation}
Another interesting consequence of Corollary 1 is the following.
If we define by continuity $f(0)=\e \eta = \e \Phi(h+\sigma z)$ and write
$\lambda = \lambda\cdot 1 + (1-\lambda)\cdot 0$ then convexity of $f(m)$ 
implies
\begin{equation}
\e \exp(\lambda \eta) \leq (\e \exp \eta )^{\lambda^2} 
\exp(\lambda(1-\lambda) \e\eta).
\label{cons}
\end{equation}
If $A=\log \e \exp(\eta -\e \eta)<\infty$
then Chebyshev's inequality and (\ref{cons}) imply that
$$
\p(\eta\geq \e\eta + t) \leq 
\e \exp(\lambda \eta -\lambda \e\eta -\lambda t)
\leq \exp(\lambda^2 A - \lambda t)
$$
and minimizing over $\lambda\in [0,1]$ we get,
\begin{equation}
\p(\eta\geq \e\eta + t) \leq
\left\{
\begin{array}{c c}
\exp(-t^2/4A), & t\leq 2A \\ 
\exp(A-t), & t\geq 2A. 
\end{array}
\right.
\label{concentr}
\end{equation}
This result can be slightly generalized. 
\begin{corollary}
If $\eta = \Phi(|\vec{h} +\vec{z}|)$ for some $\vec{h}\in\Reals^n$ 
and standard Gaussian $\vec{z}\in \Reals^n$ then the function
$m^{-1}\log \e \exp m\eta$ is convex in $m$ and, thus,
(\ref{cons}) and (\ref{concentr}) hold.
\end{corollary}

The proof follows along the lines of the proof of Corollary 1 
(or Theorem \ref{T2} in the simplest case of Corollary 1) 
and will be omitted.

\section{Proof of Theorem \ref{T2}.}

The proof of Theorem \ref{T2} will be based on the following 
observations.
First of all, we will compute the derivative of 
$\P_k$ with respect to $q_l.$ We will need the following notations.
For $0\leq l\leq k$ we define 
\begin{equation}
V_l = V_l(x,z_l) = \exp m_l (\Phi_{l+1}(x+z_l)-\Phi_l(x)).
\label{Vl}
\end{equation}
Let $Z=h+z_0+\ldots+z_k$ and $Z_l = h+z_0+\ldots +z_{l-1}$
and define
$$
X_l=\Phi_l(Z_l) \mbox{ and } W_l = V_l(Z_l,z_l) = \exp m_l(X_{l+1} -X_l).
$$
Then  the following holds.
\begin{lemma}\label{L1}
For $1\leq l\leq k,$ we have,
\begin{equation}
\frac{\partial \P_k}{\partial q_l} = 
- \frac{1}{2}(m_l-m_{l-1})\xi''(q_l) U_l
\label{Dql}
\end{equation}
where
\begin{equation}
U_l = U_l(\vec{m}, \vec{q})
= \e W_1\ldots W_{l-1}\Bigl(\e_l W_l\ldots W_k \Phi'(Z)\Bigr)^2.
\label{Ul}
\end{equation}
\end{lemma}
{\bf Proof.}
The proof of Lemma is rather straightforward (it is based
on Gaussian integration by parts) and will be omitted (see,
for example, \cite{SG} or \cite{Guerra}). 
\qed

It turns out that the function $U_l$ is nondecreasing in each $m_j$
which is the main ingredient in the proof of Theorem \ref{T2}.

\begin{theorem}\label{T1}
For any $1\leq l\leq k$ the function $U_l$ defined in (\ref{Ul})
is nondecreasing in each $m_j$ for $1\leq j\leq k.$
\end{theorem}

First, let us show how Lemma \ref{L1} and Theorem \ref{T1}
imply Theorem \ref{T2}.

{\bf Proof of  Theorem \ref{T2}.}
Let us assume that $n_j \leq m_j$ 
for all $j\leq k.$ The opposite case can be handled similarly.
If we define
$$
\vec{m}^l = (n_0,\ldots,n_l,m_{l+1},\ldots,m_k)
$$
then
$$
\P_k(\vec{n}) - \P_k(\vec{m}) = \sum_{0\leq l\leq k}
\bigl(\P_k(\vec{m}^l) - \P_k(\vec{m}^{l-1})\bigr).
$$
We will prove that
\begin{equation}
\P_k(\vec{m}^l) - \P_k(\vec{m}^{l-1}) \geq
\frac{\partial \P_k(\vec{m})}{\partial m_l}(n_l - m_l)
\label{coord}
\end{equation}
which, obviously, will prove Theorem \ref{T2}.
Let us consider vectors
$$
\vec{m}_{+}^l = (n_0,\ldots,n_{l}, m_l, m_{l+1},\ldots,m_k)
$$
and
$$
\vec{q}^l(t)=(q_0,\ldots, q_l,q_{l+1}(t), q_{l+1},q_{l+2},\ldots,q_k),
$$
where $q_{l+1}(t) = q_l + t(q_{l+1} - q_l).$
Notice that we inserted one coordinate in vectors 
$\vec{m}^l$ and $\vec{q}.$ For $0\leq t\leq 1,$ we consider
$$
\varphi(t) = \P_{k+1}(\vec{m}_{+}^l,\vec{q}^l(t)).
$$
It is easy to see that $\varphi(t)$ interpolates between
$\varphi(1) = \P_k(\vec{m}^l)$ and 
$\varphi(0) = \P_k(\vec{m}^{l-1}).$ By Lemma \ref{L1},
$$
\varphi'(t) =
%\frac{\partial \P_{k+1}(t)}{\partial t} = 
-\frac{1}{2}(m_l-n_{l})\xi''(q_{l+1}(t)) U_{l+1}
$$
where $U_{l+1}$ is defined in terms of $\vec{m}_{+}^l$
and $\vec{q}^l(t).$
Next, let us consider
$$
\vec{m}^l_{\eps} = (m_0,\ldots,m_{l-1}, m_l -\eps(m_l -n_l),
m_l,m_{l+1},\ldots,m_k)
$$
and define
$$
\varphi_{\eps}(t) = \P_{k+1}(\vec{m}_{\eps}^l,\vec{q}^l(t)).
$$
First of all, we have $\varphi_{\eps}(0) = \P_k(\vec{m})$ and 
$\varphi_{\eps}(1) = \P_k(\vec{m}_{\eps}),$ where
$$
\vec{m}_{\eps} = (m_0,\ldots,m_{l-1}, m_l -\eps(m_l-n_l),m_{l+1},\ldots,m_k).
$$
Again, by Lemma \ref{L1},
$$
\varphi_{\eps}'(t) =
%\frac{\partial \P_{k+1}^{\eps}(t)}{\partial t} = 
-\frac{1}{2}\eps(m_l-n_{l})\xi''(q_{l+1}(t)) U_{l+1}^{\eps}
$$
where $U_{l+1}^{\eps}$ is defined in terms of $\vec{m}_{\eps}^l$
and $\vec{q}^l(t).$
It is obvious that for $\eps\in[0,1]$ each coordinate
of $\vec{m}^l_{\eps}$ is not smaller than the corresponding
coordinate of $\vec{m}^l$ and, therefore, Theorem \ref{T1}
implies that $U_{l+1} \leq U_{l+1}^{\eps}.$ This implies
$$
\frac{1}{\eps}\varphi_{\eps}'(t)
\leq
\varphi'(t)
%\frac{1}{\eps}\frac{\partial \P_{k+1}^{\eps}(t)}{\partial t}
%\leq
%\frac{\partial \P_{k+1}(t)}{\partial t}
$$
and, therefore,
$$
\frac{1}{\eps}(\varphi_{\eps}(1) - \varphi_{\eps}(0))\leq
\varphi(1) - \varphi(0)
$$
which is the same as
$$
\frac{1}{\eps}(\P_k(\vec{m}_{\eps}) - \P_k(\vec{m})) \leq
\P_k(\vec{m}^l) - \P_k(\vec{m}^{l-1}).
$$
Letting $\eps\to 0$ implies (\ref{coord}) and this finishes
the proof of Theorem \ref{T2}.

\qed

\section{Proof of Theorem \ref{T1}.}

Let us start by proving some preliminary results.
Consider two classes of (smooth enough) functions 
\begin{equation}
\C=\{f:\Reals\to[0,\infty) : f(-x)=f(x), f'(x)\geq 0 \mbox{ for } x\geq 0\}
\label{C}
\end{equation}
and
\begin{equation}
\Cp=\{f:\Reals\to[0,\infty) : f(-x)=-f(x), f'(x)\geq 0 \mbox{ for } x\geq 0\}.
\label{Cprime}
\end{equation}
The next Lemma describes several facts that will be
useful in the proof of Theorem \ref{T1}.

\begin{lemma}\label{L2} 
For all $1\leq l\leq k$ and $V_l = V_l(x,z_l)$
defined in (\ref{Vl}) we have,

(a) $\Phi_l(x)$ is convex, $\Phi_l(x) \in \C$ and 
$$
\Phi_l'(x) = \e_l V_l \ldots V_k \Phi'(x+z_l+\ldots+z_k)
\in \Cp.
$$

(b) If $f_1\in\C$ and $f_2\in\Cp$ then for $x\geq 0$
$$
\e_l V_l f_1(x+z_l)f_2(x+z_l) \geq \e_l V_l f_1(x+z_l)
\e_l V_l f_2(x+z_l).
$$

(c) If $f(-x)=-f(x)$ and $f(x)\geq 0$ for $x\geq 0$
then $g(x)= \e_l V_l f(x+z_l)$ is such that
$$
g(-x) = -g(x) \mbox{ and } g(x)\geq 0 \mbox{ if } x\geq 0.
$$ 

(d) If $f\in\C$ then $\e_l V_l f(x+z_l)\in \C.$
\hspace{0.2cm} (e) If $f\in\Cp$ then $\e_l V_l f(x+z_l)\in \Cp.$

(f) $f(x)=\e_l V_l\log V_l \in \C. $

\end{lemma}

{\bf Proof.}
(a) Since $\Phi_{k+1}$ is convex, symmetric and nonnegative then 
$\Phi_l(x)$ is convex, symmetric and nonnegative by induction on $l$
in (\ref{rec}). Convexity is the consequence of
H\"older's inequality and symmetry follows from
\begin{eqnarray*}
\Phi_l(-x) = \frac{1}{m_l}\log \e_l \exp m_l \Phi_{l+1}(-x+z_l)
&=& 
\frac{1}{m_l}\log \e_l \exp m_l \Phi_{l+1}(-x-z_l)
\\
&=& 
\frac{1}{m_l}\log \e_l \exp m_l \Phi_{l+1}(x+z_l)
= \Phi_l(x).
\end{eqnarray*}
Obviously, this implies that $\Phi_l'(x) \in \Cp.$ 

(b) Let $z_l'$ be an independent copy of $z_l$ and, for simplicity
of notations, let $\sigma^2=\e z_l^2$. Since $\e_l V_l = 1$ 
(i.e. we can think of $V_l$ as the change of density),
we can write,
\begin{eqnarray}
&&
\e_l V_l f_1(x+z_l) f_2(x+z_l) - \e_l V_l f_1(x+z_l) \e_l V_l f_2(x+z_l)=
\label{b1}
\\
&=&
\e_l V_l(x,z_l) V_l(x,z_l')\Bigl(f_1(x+z_l) - f_1(x+z_l')\Bigr)
\Bigl(f_2(x+z_l)-f_2(x+z_l')\Bigr)
I(z_l\geq z_l')
\nonumber
\end{eqnarray}
Since 
$V_l(x,z_l)V_l(x,z_l') = \exp m_l(\Phi_l(x+z_l) + \Phi_l(x+z_l') -2\Phi_l(x)),$
if we make the change of variables $s=x+z_l$ and $t=x+z_l'$ then
the right hand side of (\ref{b1}) can be written as
\begin{eqnarray}
\frac{1}{2\pi \sigma^2} \exp (-2m_l\Phi_l(x))
\int\limits_{\{s\geq t\}} K(s,t)
\exp\Bigl(-\frac{1}{2\sigma^2}((s-x)^2 +(t-x)^2)\Bigr)
ds dt,
\label{Kint}
\end{eqnarray}
where
$$
K(s,t) = \exp m_l(\Phi_l(s) + \Phi_l(t))
\Bigl(f_1(s) - f_1(t)\Bigr)\Bigl(f_2(s)-f_2(t)\Bigr).
$$
We will split the region of integration 
$\{s\geq t\}=\Omega_1\cup \Omega_2$
in the last integral into two disjoint sets 
$$
\Omega_1 = \{(s,t) : s\geq t, |s|\geq |t|\},\,\,\,
\Omega_2 = \{(s,t) : s\geq t, |s| < |t|\}.
$$
In the integral over $\Omega_2$ we will make the change
of variables $s=-v, t=-u$ so that for $(s,t)\in \Omega_2$
we have $(u,v)\in \Omega_1$ and $ds dt = du dv.$
Also,
$$
K(s,t) = K(-v,-u) = -K(u,v)
$$
since $\Phi_l$ is symmetric by (a), $f_1\in \C, f_2\in\Cp$
and, therefore, 
$$
\Bigl(f_1(-v)-f_1(-u)\Bigr)\Bigl(f_2(-v) - f_2(-u)\Bigr)
=
-\Bigl(f_1(u)-f_1(v)\Bigr)\Bigl(f_2(u) - f_2(v)\Bigr).
$$
Therefore,
$$
\int\limits_{\Omega_2} K(s,t)
\exp\Bigl(-\frac{1}{2\sigma^2}((s-x)^2 +(t-x)^2)\Bigr)
ds dt
=
-\int\limits_{\Omega_1} K(u,v)
\exp\Bigl(-\frac{1}{2\sigma^2}((u+x)^2 +(v+x)^2)\Bigr)
du dv
$$
and (\ref{Kint}) can be rewritten as
\begin{equation}
\frac{1}{2\pi \sigma^2} \exp (-2m_l\Phi_l(x))
\int\limits_{\Omega_1} K(s,t) L(s,t,x) ds dt
\label{Kint2}
\end{equation}
where
$$
L(s,t,x) = \exp\Bigl(-\frac{1}{2\sigma^2}((s-x)^2 +(t-x)^2)\Bigr)
-\exp\Bigl(-\frac{1}{2\sigma^2}((s+x)^2 +(t+x)^2)\Bigr).
$$
Since $f_1\in\C,$ for $(s,t)\in\Omega_1$ we have $f_1(s)-f_1(t)
= f_1(|s|) - f_1(|t|) \geq 0.$ Moreover, since for $(s,t)\in\Omega_1$
we have $t\leq s,$ the fact that $f_2\in\Cp$
implies that $f_2(s) - f_2(t) \geq 0.$ Combining these two
observations we get that $K(s,t) \geq 0$ on $\Omega_1.$
Finally, for $(s,t)\in\Omega_1$ we have $L(s,t,x)\geq 0$
because 
$$
(s-x)^2 + (t-x)^2 \leq (s+x)^2 + (t+x)^2
\Longleftrightarrow
x(s+t)\geq 0,
$$
and the latter holds because $x\geq 0$ and
$s+t\geq 0$ on $\Omega_1.$
This proves that (\ref{Kint2}), (\ref{Kint})
and, therefore, the right hand side of (\ref{b1}) are
nonnegative. 

(c) Let $g(x) = \e_l V_l(x,z_l) f(x+z_l).$ Then
$$
g(-x) = \e_l V_l(-x,z_l) f(-x+z_l) = \e_l V(-x,-z_l)
f(-x-z_l) = -\e_l V_l(x,z_l) f(x+z_l) = -g(x). 
$$
Next, if $x\geq 0$ and $\sigma^2=\e_l z_l^2$ then
\begin{eqnarray*}
g(x)
&=& 
\exp(-m_l \Phi_l'(x)) \e_l \exp( m_l \Phi_{l+1}(x+z_l)) f(x+z_l)
= \exp(-m_l \Phi_l'(x)) \frac{1}{\sqrt{2\pi} \sigma} \times
\\
&& 
\times
\int_{s\geq 0} \exp (m_l \Phi_{l+1}(s)) f(s) \Bigl(
\exp\bigl( -\frac{1}{2\sigma^2} (x-s)^2\bigr) 
- \exp\bigl( -\frac{1}{2\sigma^2} (x+s)^2\bigr) 
\Bigr) ds \geq 0
\end{eqnarray*}
because $(x-s)^2 \leq (x+s)^2$  for $x,s\geq 0$
and $f(s)\geq 0$ for $s\geq 0.$

(d) Take $f\in \C.$ Positivity of $\e_l V_l f(x+z_l)$
is obvious and symmetry follows from
\begin{equation}
\e_l V_l(-x,z_l)f(-x+z_l) = \e_l V_l(-x,-z_l) f(-x-z_l)
=\e_l V_l(x,z_l) f(x+z_l).
\label{symex}
\end{equation}
Let $x\geq 0.$
Recalling the definition (\ref{Vl}), the derivative 
$$
\frac{\partial}{\partial x} \e_l V_l(x,z_l) f(x+z_l)
= \I+ m_l \II
$$
where $\I = \e_l V_l(x,z_l) f'(x+z_l)$ and 
\begin{eqnarray*}
\II 
&=& 
\e_l V_l(x,z_l) f(x+z_l)
(\Phi_{l+1}'(x+z_l) - \Phi_l'(x))
\\
&=&
\e_l V_l(x,z_l) f(x+z_l) \Phi_{l+1}'(x+z_l)
- \e_l V_l(x,z_l) f(x+z_l) \e_l V_l \Phi_{l+1}'(x+z_l),
\end{eqnarray*}
since (\ref{rec}) yields that 
$\Phi_l'(x) = \e_l V_l(x,z_l) \Phi_{l+1}'(x+z_l).$
By (a), $\Phi_{l+1}'\in\Cp,$ and since $f\in\C,$ (b)
implies that $\II\geq 0.$
The fact that $\I \geq 0$ for $x\geq 0$
follows from (c) because $f'(-x) = -f'(x)$ and
$f'(x)\geq 0$ for $x\geq 0.$

(e) Take $f\in \Cp.$ Antisymmetry of $\e_l V_l f(x+z_l)$
follows from
$$
\e_l V_l(-x,z_l)f(-x+z_l) = \e_l V_l(-x,-z_l) f(-x-z_l)
= - \e_l V_l(x,z_l) f(x+z_l).
$$
As in (d), the derivative can be written as
$$
\frac{\partial}{\partial x} \e_l V_l(x,z_l) f(x+z_l)
= \I+ m_l \II
$$
where $\I = \e_l V_l(x,z_l) f'(x+z_l)$ and
$$
\II =
\e_l V_l(x,z_l) f(x+z_l) \Phi_{l+1}'(x+z_l)
- \e_l V_l(x,z_l) f(x+z_l) \e_l V_l \Phi_{l+1}'(x+z_l).
$$
First of all, $\I \geq 0$ because $f'\geq 0$ for $f\in\Cp.$
As in (\ref{b1}) we can write
$$
\II = \e_l V_l(x,z_l) V_l(x,z_l')
\bigl(f(x+z_l) - f(x+z_l')\bigr)
\bigl(\Phi_{l+1}'(x+z_l)- \Phi_{l+1}'(x+z_l')\bigr)
I(z_l\geq z_l').
$$
But both $f$ and $\Phi_{l+1}'$ are in the class $\Cp$
and, therefore, both nondecreasing
which, obviously, implies that they are similarly ordered,
i.e. for all $a,b\in \Reals,$
\begin{equation}
(f(a) - f(b))(\Phi_{l+1}'(a)- \Phi_{l+1}'(b)) \geq 0
\label{sim}
\end{equation}
and as a result $\II \geq 0.$

(f) Symmetry of $g(x)= \e_l V_l \log V_l$ follows as above
and positivity follows from Jensen's inequality, convexity
of $x\log x$ and the fact that $\e_l V_l = 1.$
Next, using that $\Phi_l'(x) = \e_l V_l \Phi_{l+1}'(x+z_l)$
we can write
\begin{eqnarray*}
g'(x) 
&=& 
m_l \e_l (1+\log V_l) V_l (\Phi_{l+1}'(x+z_l) - \Phi_l'(x))
\\
&=& 
m_l^2 \e_l V_l (\Phi_{l+1}(x+z_l) - \Phi_l(x))
(\Phi_{l+1}'(x+z_l) - \Phi_l'(x))
\\
&=& 
m_l^2 \e_l V_l \Phi_{l+1}(x+z_l) 
(\Phi_{l+1}'(x+z_l) - \Phi_l'(x))
\\
&=&
m_l^2 \Bigl(
\e_l V_l \Phi_{l+1}(x+z_l) \Phi_{l+1}'(x+z_l) 
- \e_l V_l \Phi_{l+1}(x+z_l) \e_l V_l \Phi_{l+1}'(x+z_l) 
\Bigr).
\end{eqnarray*}
Since $\Phi_{l+1}\in\C$ and $\Phi_{l+1}'\in\Cp,$
(b) implies that for $x\geq 0,$ $g'(x)\geq 0$
and, therefore, $g\in\C.$

\qed

{\bf Proof of Theorem \ref{T1}.}

We will consider two separate cases.

{\bf Case 1}. $j\leq l-1.$ 
First of all, using Lemma \ref{L2} (a) we can rewrite
$U_l$ as
$$
U_l = \e W_1\ldots W_{l-1} f_l(Z_l)
$$
where 
\begin{equation}
f_l(x) = (\Phi_l'(x))^2 \in \C \mbox{ since }
\Phi_l'(x)\in \Cp.
\label{fell}
\end{equation}
Using that
$$
X_j = \frac{1}{m_j}\log \e_j \exp m_j X_{j+1}
$$
we get
$$
\frac{\partial X_j}{\partial m_j} 
= 
\frac{1}{m_j}\e_j W_j X_{j+1} - \frac{1}{m_j^2} \log \e_j \exp m_j X_{j+1}
= \frac{1}{m_j}\e_j W_j(X_{j+1} - X_j). 
$$
For $p\leq j,$ we get
$$
\frac{\partial X_p}{\partial m_j} = 
\frac{1}{m_j}\e_p W_p \ldots W_j (X_{j+1}-X_j), 
$$
and for $p> j,$ $X_p$ does not depend on $m_j.$
Therefore,
\begin{eqnarray*}
&&
\frac{\partial}{\partial m_j} W_1 \ldots W_{l-1} 
=
\frac{\partial}{\partial m_j}
\exp \Bigl( \sum_{p\leq l-1} m_p(X_{p+1} -X_p) \Bigr)
\\
&&
=
W_1\ldots W_{l-1}\Bigl(
(X_{j+1} - X_j) - \frac{1}{m_j}\sum_{p\leq j}(m_p - m_{p-1})
\e_p W_p \ldots W_j (X_{j+1} - X_j) \Bigr).
\end{eqnarray*}
Hence,
\begin{eqnarray*}
m_j \frac{\partial U_l}{\partial m_j}
&=&
m_j \e W_1\ldots W_{l-1} f_l(Z_l)(X_{j+1} - X_j) 
\\
&&
-\sum_{p\leq j}(m_p - m_{p-1}) \e W_1 \ldots W_{l-1}
f_l (Z_l) \e_p W_p \ldots W_{j}(X_{j+1} - X_j).
\end{eqnarray*}
If we denote $f_j(Z_{j+1}) = \e_{j+1} W_{j+1}\ldots W_{l-1} f_l(Z_l)$
then we can rewrite
\begin{eqnarray}
m_j \frac{\partial U_l}{\partial m_j}
&=&
m_j \e W_1\ldots W_{j} f_j(Z_{j+1})(X_{j+1} - X_j) 
\label{case1}
\\
&&
-\sum_{p\leq j}(m_p - m_{p-1}) \e W_1 \ldots W_{p-1}
\e_p W_p \ldots W_j f_j(Z_{j+1}) \e_p W_p \ldots W_{j}(X_{j+1} - X_j).
\nonumber
\end{eqnarray}
First of all, let us show that
\begin{equation}
\e_j W_j f_j(Z_{j+1})(X_{j+1} - X_j) \geq
\e_j W_j f_j(Z_{j+1}) \e_j W_j(X_{j+1} - X_j).
\label{ina0}
\end{equation}
Since $X_j$ does not depend on $z_j$ and $\e_j W_j = 1,$ 
this is equivalent to
\begin{equation}
\e_j W_j f_j(Z_{j+1}) X_{j+1}  \geq
\e_j W_j f_j(Z_{j+1}) \e_j W_j X_{j+1}.
\label{ina}
\end{equation}
Here $f_{j}$ and $X_{j+1}$ are both functions of 
$Z_{j+1} = Z_j +z_j.$  Since by (\ref{fell}), 
$f_l(Z_l)$ seen as a function of $Z_l$ is in $\C,$
applying Lemma \ref{L2} (d) inductively we get
that $f_j(Z_{j+1})$ seen as a function of $Z_{j+1}$
is also in $\C.$ By Lemma \ref{L2} (a), $X_{j+1}$
seen as a function of $Z_{j+1}$ is also in $\C.$
Therefore, $f_j$ and $X_{j+1}$ are similarly ordered
i.e.
$$
\bigl(f_j(Z_{j+1}) - f_j(Z_{j+1}')\bigr)
\bigl(X_{j+1}(Z_{j+1}) - X_{j+1}(Z_{j+1}') \bigr) \geq 0
$$
and, therefore, using the
same trick as in (\ref{b1}) (often called Chebyshev's
inequality) we get (\ref{ina}) and, hence, (\ref{ina0}). 
By Lemma \ref{L2} (d), $\e_j W_j f_j(Z_{j+1})$ seen as a function
of $Z_j$ is in $\C$ and by Lemma \ref{L2} (f),
$\e_j W_j (X_{j+1} - X_j) = m_j^{-1} \e_j W_j \log W_j$
seen as a function of $Z_j$ is also in $\C.$
Therefore, they are similarly ordered and again 
$$
\e_p W_p \ldots W_{j-1} \e_j W_j f_j(Z_{j+1}) \e_j W_j (X_{j+1} - X_j) 
\geq
\e_p W_p\ldots W_j f_j(Z_{j+1}) \e_pW_p\ldots W_j(X_{j+1} - X_j).
$$
Combining this with (\ref{ina0}) implies that
$$
\e W_1\ldots W_{j} f_j(Z_{j+1})(X_{j+1} - X_j) \geq
\e W_1 \ldots W_{p-1}
\e_p W_p \ldots W_j f_j(Z_{j+1}) 
\e_p W_p \ldots W_{j}(X_{j+1} - X_j).
$$
Since $m_j = \sum_{p\leq j}(m_p - m_{p-1}),$ this and
(\ref{case1}) imply that $\partial U_l/\partial m_j \geq 0$
which completes the proof of Case 1. 
\qed

{\bf Case 2.} $j\geq l.$ If we denote 
$$
g_l = g_l(Z_l) = \e_l W_l \ldots W_l \Phi'(Z),\,\,\,
f_l = f_l(Z_l) = g_l^2
$$
then a straightforward calculation similar to the one
leading to (\ref{case1}) gives
\begin{eqnarray}
m_j \frac{\partial U_l}{\partial m_j} 
=
&-&
\sum_{p\leq l-1} (m_{p}-m_{p-1}) \e W_1 \ldots W_{l-1} f_l 
\e_p W_p \ldots W_j (X_{j+1} - X_j)
\nonumber
\\
&-&
(2m_l - m_{l-1}) \e W_1 \ldots W_{l-1} f_l 
\e_l W_l \ldots W_j (X_{j+1} - X_j)
\nonumber
\\
&-&
\sum_{l+1\leq p\leq j} 2(m_{p} - m_{p-1}) \e W_1\ldots W_{l-1} g_l
\e_l W_l \ldots W_k \Phi'(Z) \e_p W_p\ldots W_j (X_{j+1} -X_j)
\nonumber
\\
&+&
2m_j \e W_1\ldots W_{l-1} g_l
\e_l W_l \ldots W_k \Phi'(Z) (X_{j+1} -X_j).
\label{case2}
\end{eqnarray}
To show that this is positive we notice that
$$
2m_j = \sum_{p\leq l-1} (m_{p}-m_{p-1}) + (2m_l - m_{l-1})
+\sum_{l+1\leq p\leq j} 2(m_{p} - m_{p-1})
$$
and we will show that the last term with factor $2m_j$
is bigger than all other terms with negative factors.
If we denote 
$$
h(Z_{j+1}) = \e_{j+1} W_{j+1} \ldots W_k \Phi'(Z) 
$$
then since $\Phi' \in \Cp,$ using Lemma \ref{L2} (e)
inductively, we get that $h(Z_{j+1})$ seen as a function
of $Z_{j+1}$ is in $\Cp.$ Each term in the third line of
(\ref{case2}) (without the factor $2(m_p - m_{p-1})$)
can be rewritten as
\begin{equation}
\e W_1\ldots W_{l-1} g_l \e_l W_l \ldots W_{p-1} 
\e_p W_p \ldots W_j h(Z_{j+1}) \e_p W_p\ldots W_j (X_{j+1} -X_j),
\label{comp1}
\end{equation}
the term in the second line of (\ref{case2}) 
(without the factor $2m_l - m_{l-1})$) is equal to (\ref{comp1})
for $p=l,$
and the term in the fourth line (without $2m_j$) can be written as 
\begin{equation}
\e W_1\ldots W_{l-1} g_l \e_l W_l \ldots W_j h(Z_{j+1}) (X_{j+1} -X_j).
\label{comp2}
\end{equation}
We will show that (\ref{comp2}) is bigger than (\ref{comp1})
for $l\leq p\leq j.$ This is rather straightforward using Lemma \ref{L2}.
Notice that $g_l = g_l(Z_l)$ seen as a function of $Z_l$
is in $\Cp$ by Lemma \ref{L2} (a). If we define for $l\leq p \leq j,$
$$
r_p(Z_l) = \e_l W_l \ldots W_{p-1} 
\e_p W_p \ldots W_j h(Z_{j+1}) \e_p W_p\ldots W_j (X_{j+1} -X_j)
$$
and 
$$
r(Z_l) = \e_l W_l \ldots W_j h(Z_{j+1}) (X_{j+1} -X_j)
$$
then the difference of (\ref{comp2}) and (\ref{comp1}) is
\begin{equation}
\e W_1\ldots W_{l-1} g_l(Z_l)(r(Z_l) - r_p(Z_l)).
\label{diff}
\end{equation}
Using the argument similar to (\ref{symex}) (and several
other places above), it should be obvious that 
$r_p(-Z_l)=-r_p(Z_l)$ since $X_i$'s are symmetric and $h$
is antisymmetric. Similarly, $r(-Z_l) = -r(Z_l).$
Therefore, if we can show that 
\begin{equation}
r(Z_l) - r_p(Z_l)\geq 0 \mbox{ for } Z_l\geq 0
\label{diff2}
\end{equation}
then, since $g_l\in\Cp,$ we would get that
$$
g_l(Z_l)(r(Z_l) - r_p(Z_l)) \geq 0 \mbox{ for all } Z_l
$$
and this would prove that (\ref{diff}) is nonnegative.
Let us first show that (\ref{diff2}) holds for $p=j.$ 
In this case, since $X_j$ does not depend on $z_j$
and, therefore, $\e_j W_j X_j = X_j,$ (\ref{diff2})
is equivalent to
\begin{equation}
\e_l W_l \ldots W_{j-1} \e_j W_j h(Z_{j+1}) X_{j+1} 
\geq 
\e_l W_l \ldots W_{j-1} \e_j W_j h(Z_{j+1}) \e_j W_j X_{j+1},
\label{diff3}
\end{equation}
for $Z_l \geq 0.$
Let us define
$$
\Delta_j(Z_j) = \e_j W_j h(Z_{j+1})X_{j+1} - 
\e_j W_j h(Z_{j+1}) \e_j W_j X_{j+1}.
$$
As above, $\Delta_j(-Z_j) = -\Delta_j(Z_j)$
and by Lemma \ref{L2} (b), $\Delta_j(Z_j) \geq 0$
for $Z_j \geq 0,$ since  $h\in\Cp$ and $X_{j+1}\in \C.$
Therefore, by Lemma \ref{L2} (c),
$$
\Delta_{j-1}(Z_{j-1}):= \e_{j-1} W_{j-1} \Delta_{j}(Z_{j-1} +z_j)
\geq 0 \mbox{ if } Z_{j-1}\geq 0
$$
and, easily, $\Delta_{j-1}(-Z_{j-1}) = -\Delta_{j-1}(Z_{j-1}).$
Therefore, if for $i\geq l$ we define 
$$
\Delta_i(Z_i) = \e_i W_i \Delta_{i+1}(Z_i+z_i)
$$
we can proceed by induction to show that 
$\Delta_i(-Z_i)=-\Delta_i(Z_i)$ and $\Delta_i(Z_i)\geq 0$
for $Z_i\geq 0.$ For $i=l$ this proves (\ref{diff3}) and,
therefore, (\ref{diff2}) for $p=j.$
Next, we will show that
\begin{equation}
r_{p+1}(Z_l) - r_p(Z_l)\geq 0 \mbox{ for } Z_l\geq 0
\label{diff4}
\end{equation}
for all $l\leq p<j,$ and this, of course, will prove (\ref{diff2}).
If we define 
$$
f_1(Z_{p+1})=\e_{p+1} W_{p+1}\ldots W_j h(Z_{j+1}) \mbox{ and }
f_2(Z_{p+1})=\e_{p+1} W_{p+1}\ldots W_j (X_{j+1} - X_j)
$$
then (\ref{diff4}) can be rewritten as
$$
\e_l W_l \ldots W_{p-1} \e_p W_p f_1(Z_{p+1}) f_2(Z_{p+1})
\geq 
\e_l W_l \ldots W_{p-1} \e_p W_p f_1(Z_{p+1}) \e_p f_2(Z_{p+1})
\mbox{ for }  Z_l\geq 0. 
$$
Since $h(Z_{j+1})\in \Cp,$ recursive application of Lemma \ref{L2} (e)
implies that $f_1(Z_{p+1})\in\Cp.$ Since 
$\e_j W_j (X_{j+1} -X_j)= m_j^{-1} \e_j W_j \log W_j$ seen
as a function of $Z_j$ is in $\C$ by Lemma \ref{L2} (f),
recursive application of Lemma \ref{L2} (d) implies that
$f_2(Z_{p+1}) \in \C.$ 
If we now define
$$
\Delta_p(Z_p) = \e_p W_p f_1(Z_{p+1}) f_2(Z_{p+1}) - 
\e_p W_p f_1(Z_{p+1}) \e_p W_p f_2(Z_{p+1}),
$$
then, as above, $\Delta_p(-Z_p) = -\Delta_p(Z_p)$
and by Lemma \ref{L2} (b), $\Delta_p(Z_p) \geq 0$
for $Z_p \geq 0,$ since  $f_1\in\Cp$ and $f_2\in \C.$
Therefore, by Lemma \ref{L2} (c),
$$
\Delta_{p-1}(Z_{p-1}):= \e_{p-1} W_{p-1} \Delta_{p}(Z_{p-1} +p_j)
\geq 0 \mbox{ if } Z_{p-1}\geq 0
$$
and, easily, $\Delta_{p-1}(-Z_{p-1}) = -\Delta_{p-1}(Z_{p-1}).$
Therefore, if for $i\geq l$ we define 
$$
\Delta_i(Z_i) = \e_i W_i \Delta_{i+1}(Z_i+z_i)
$$
we can proceed by induction to show that 
$\Delta_i(-Z_i)=-\Delta_i(Z_i)$ and $\Delta_i(Z_i)\geq 0$
for $Z_i\geq 0.$ For $i=l$ this proves (\ref{diff4}). 
Thus, we finally proved that (\ref{comp2}) is bigger
than (\ref{comp1}) for $p\geq l.$ To prove that (\ref{case2})
is nonnegative it remains to show that each term in the first line
of (\ref{case2}) (with the factor $m_p - m_{p-1}$) is smaller
than (\ref{comp2}) it is enough to show that
\begin{equation}
\e W_1 \ldots W_{l-1} f_l \e_p W_p \ldots W_j (X_{j+1} - X_j)
\leq
\e W_1 \ldots W_{l-1} f_l \e_l W_l \ldots W_j (X_{j+1} - X_j)
\label{la}
\end{equation}
since the right hand side is (\ref{comp1}) for $p=l$
and it was already shown to be smaller than (\ref{comp2}).
But the proof of (\ref{la}) follows exactly the same argument
as the proof of (\ref{ina0}) in Case 1 and this finishes
the proof of Case 2.

\qed


\begin{thebibliography}{99} 


\bibitem{Guerra} Guerra, F. (2003) Broken replica 
symmetry bounds in the mean field spin glass model. 
{\it Comm. Math. Phys.} {\bf 233}, no. 1, 1-12.


\bibitem{SherK} Sherrington, D., Kirkpatrick, S. (1972)
Solvable model of a spin glass. 
{\it Phys. Rev. Lett.} {\bf 35}, 1792-1796.

\bibitem{SG} Talagrand, M. (2003)   
Spin Glasses: a Challenge for Mathematicians. 
Springer-Verlag.


\bibitem{T-gP} Talagrand, M. (2003)
The generalized Parisi formula. 
{\it C. R. Math. Acad. Sci. Paris} {\bf 337}, no. 2, 111-114.

\bibitem{T-P} Talagrand, M. (2003)
Parisi formula. To appear in {\it Ann. Math.}

\bibitem{T-M} Talagrand, M. (2003)
On the meaning of Parisi's functional order parameter.
{\it C. R. Math. Acad. Sci. Paris} {\bf 337}, no. 9, 625-628.

\bibitem{T-Par} Talagrand, M. (2004)
Parisi measures. Preprint.



\end{thebibliography}
\end{document}